\documentclass[12pt]{article}
\usepackage[utf8]{inputenc}
\usepackage[english]{babel}
\usepackage{color}
\usepackage{graphicx}
\usepackage{vmumath}
\usepackage{amsmath}
\usepackage{amsthm}
\usepackage{avo_book}
\usepackage{hyperref}
\usepackage{caption}
\usepackage{booktabs}
\usepackage{subcaption}
\usepackage[a4paper,left=15mm,right=15mm,top=30mm,bottom=20mm]{geometry}
\noaftermath
\begin{document}

\bibliographystyle{unsrt}

\title{Enumeration of $r$-regular Maps on the Torus. \\Part I: Enumeration of Rooted and Sensed Maps}

\author{
Evgeniy Krasko \qquad  Alexander Omelchenko\\
\small St. Petersburg Academic University\\
\small 8/3 Khlopina Street, St. Petersburg, 194021, Russia\\
\small\tt \{krasko.evgeniy, avo.travel\}@gmail.com
}

\begin{abstract}
The work that consists of two parts is devoted to the problem of enumerating unrooted $r$-regular maps on the torus up to all its symmetries. We begin with enumerating near-$r$-regular rooted maps on the torus, projective plane and the Klein bottle. We also present the results of enumerating some special kinds of maps on the sphere: near-$r$-regular maps, maps with multiple leaves and maps with multiple root semi-edges. For $r=3$ and $r=4$ we obtain exact analytical formulas. For larger $r$ we derive recurrence relations. Then using these results we enumerate $r$-regular maps on the torus up to homeomorphisms that preserve its orientation --- so-called sensed maps. Using the concept of a quotient map on an orbifold we reduce this problem to enumeration of certain classes of rooted maps. For $r=3$ and $r=4$ we obtain closed-form expressions for the numbers of $r$-regular sensed maps by edges. All these results will be used in the second part of the work to enumerate $r$-regular maps on the torus up to all homeomorphisms --- so-called unsensed maps.  
\end{abstract}

\maketitle

\section*{Introduction}

By a (topological) map $\M$ on a surface we will mean an embedding of a connected graph $G$, loops and multiple edges allowed, into a compact $2$-dimensional surface $S$, such that $G$ is as a subset of $S$, and its complement $S\setminus G$ is homeomorphic to a set of topological polygons. These $f$ polygons are the faces of the map $\M$ which also has some amount $v$ of vertices (points on the surface $S$) and $n$ of edges (nonintersecting curves on the surface that have no common points other than the vertices of the graph). A map is rooted if one of its semi-edges is distinguished and one of two possible local orientations of the surface in the neighborhood of this semi-edge is chosen. The vertex this semi-edge is incident to will often be referred to as the root vertex. A map is called $r$-regular if the degree of each vertex is equal to $r$, and near-$r$-regular if all its vertices except possibly the root one have degree $r$. The notion of a dual map allows to consider an $r$-angulation of the surface instead of a $r$-regular map on it. 

Two topological maps $\rm M_1$ and $\rm M_2$ on a surface are said to be equivalent if there exists a homeomorphism $h$ of the surface into itself that transforms the edges, vertices and faces of one map into the corresponding elements of the other. For rooted maps we require in addition that such homeomorphism transforms the root semi-edge of $\rm M_1$ into the root semi-edge of $\rm M_2$ preserving local orientation of the surface. Two homeomorphisms $h_1$ and $h_2$ of a topological map $\rm M$ onto itself are considered equivalent if $h_1\cdot h_2^{-1}$ acts as identity on semi-edges of $\rm M$. By sensed (unsensed) map we mean an equivalence class of maps where the equivalence relation is defined by sense-preserving (sense-preserving or sense-reversing) homeomorphisms. 

William Tutte noted (see, for example \cite{Tutte_Census}) that a map with a distinguished root semi-edge always has a trivial symmetry group. In the series of his ``census" \,papers he derived formulas for enumerating rooted planar maps with $n$ edges \cite{Tutte_Census}, rooted $3$-regular maps \cite{Tutte_triangulations}, rooted eulerian maps \cite{Tutte_slicings}, and some other classes maps on the sphere. In the paper \cite{Tutte_Enum} Tutte derived a recurrence that expresses the number of rooted maps with some given numbers of vertices, faces and the degree of the root vertex through the numbers of analogous maps for smaller parameter values. He provided its solution in the form of a generating function and calculated the corresponding numbers of rooted planar maps with $v$ vertices and $f$ faces. At approximately the same time the papers \cite{Brown_disk_triangulations}, \cite{Mullin_triangular_maps_1964} and \cite{Mullin_triangular_maps_1965} appeared, devoted to the problem of enumerating $3$-regular rooted maps on the sphere. The paper \cite{Brown_quadrangulations} considered some class of $4$-regular maps. In \cite{Brown_nonplanar_graphs} Brown applied Tutte's technique to obtain the numbers of rooted non-separable maps on the sphere and rooted maps on the projective plane. 

The first work devoted to enumeration of rooted maps on surfaces of arbitrary genus $g > 0$ was the paper \cite{Walsh_Lehman} of Walsh and Lehman. Using Tutte's approach for enumerating planar maps, the authors derived a recurrence relation for the numbers of rooted maps and calculated the first terms of the corresponding sequences. In addition, they provided an explicit expression for the number of one-face (or unicellular) maps with $n$ edges on a surface of genus $g$, as well as a formula for the number of such maps with a prescribed list of vertex degrees.

The next step in enumerating rooted maps on surfaces was done by Bender and Canfield. In \cite{Bender_asymptotic_number} they derived a system of equations for generation functions that enumerates rooted maps on orientable and non-orientable surfaces, and found the asymptotics of the corresponding sequences. In \cite{Arques} D. Arques used the method of \cite{Bender_asymptotic_number} to obtain closed-form formulas for the number of rooted toroidal maps with $n$ edges and for the number of rooted toroidal maps with $v$ vertices and $f$ faces. In \cite{Bender_1988} Bender, Canfield and Robinson independently from Arques enumerated rooted maps on the torus and the projective plane, derived some explicit expressions for the corresponding generation functions and analyzed the asymptotic behavior of their coefficients. In \cite{Bender_1991} a generation function was obtained for the number of rooted maps of genera $2$ and $3$ with $n$ edges. Some recent results regarding enumeration of rooted maps on surfaces of higher genera can be found in \cite{Walsh_Giorgetti_2014} and \cite{Walsh_Giorgetti}. In 2000 appeared the Atlas \cite{An_atlas} which gives a complete list of unrooted maps and hypermaps with small number of edges on both orientable and nonorientable surfaces. For each unlabelled object the number of different rootings is also given.

Along with counting arbitrary maps, a series of works devoted to enumerating $3$- and $4$-regular rooted maps on surfaces (or their duals --- triangulations and quadrangulations) appeared in the nineties and in the beginning of the two thousandth. In 1991 Gao \cite{Gao_1991} used the approach of \cite{Bender_asymptotic_number} to enumerate rooted triangular maps on genus $g$ surfaces. In particular, he obtained a list of parametric expressions for the numbers of such maps on the sphere, torus and the projective plane. Lately Gao extended these results to some special classes of $3$-regular maps (see, for example, \cite{Gao_Wormald_2002}). Furthermore, there exists a series of recent papers \cite{Ren_Liu_Klein_bottle}, \cite{Ren_Liu_4-regular_maps_Proj_Plane}, \cite{Ren_Liu_near_4-regular_torus}, \cite{Long_Ren_2_connected} devoted to enumeration of $4$-regular maps. Unfortunately, explicit formulas obtained in these papers for arbitrary $4$-regular maps on the torus and on the Klein bottle contain some mistakes and do not align with the results of explicit generation of the corresponding structures.

A general technique for counting planar maps up to orientation-preserving homeomorphisms was developed by Liskovets \cite{Liskovets_85} in the early eighties. His approach reduces the enumerating problem for sensed maps on the sphere to counting {\em quotient maps on orbifolds}, maps on quotients of a surface under a finite group of automorphisms. His ideas were further developed by Mednykh and Nedela in a series of papers devoted to enumerating sensed maps and hypermaps on orientable surfaces (see \cite{Mednykh_Nedela}, \cite{Mednykh_Hypermaps}, \cite{Azevedo}). They employ a geometric approach based on reducing the problem for sensed and unsensed maps to enumeration of rooted maps on cyclic orbifolds. Coefficients in the formulas are determined in terms of numbers of order-preserving epimorphisms from orbifold fundamental groups to cyclic groups. 

The present paper generalizes the results of the article \cite{Krasko_Omelch_4_reg_one_face_maps} devoted to enumeration of regular rooted and sensed one-face maps on surfaces of a given genus. The first part of the present paper is devoted to building analytical formulas for enumerating near-$r$-regular rooted maps on the torus (the numbers $\tau_n$), on the projective plane ($\pi_n$) and on the Klein bottle ($\kappa_n$). For this problem we use the technique described in papers \cite{Tutte_Enum}, \cite{Walsh_Lehman} and \cite{Bender_1988}. For $r=3$ and $r=4$ we obtain simple exact analytical expressions for the corresponding numbers. For other values of $r$ we derive a system of generating functions as well as a system of recurrence relations that allows us to determine coefficients numerically. In the second part of the present article we will use the technique developed by Mednykh and Nedela to enumerate sensed $r$-regular maps on the torus. For $r=3$ and $r=4$ we provide explicit enumerating formulas. For $r=5$ and $r=6$ we give the numerical results.

\section{Counting rooted $r$-regular maps on the torus, the projective plane and the Klein bottle}

Let $t_{n,d}$ be the number of near-$r$-regular maps with $n$ edges on the torus with the root semi-edge incident to a vertex of degree $d$. To obtain recurrence relations for the numbers $t_{n,d}$ we will use the approach of Tutte \cite{Tutte_Census}, \cite{Tutte_Enum}. For the problem in question it is more convenient to use not the operation of deleting an edge, but the operation of contracting one, as it retains the near-$r$-regularity property. We show that
\begin{equation}
\label{eq:t_n_d}
t_{n,d}=t_{n-1,d+r-2}+d_{n-1,d-2}+2\sum\limits_{i=0}^{n-1}\sum\limits_{j=0}^{d-2}s_{n-i-1,d-j-2}\cdot t_{i,j},
\qquad n\geq 1,\quad d\geq 1,
\end{equation}
$$
t_{0,0}=1,\qquad\qquad t_{n,d}=0\qquad \forall\,\,n<0\quad \text{or}\quad d<0.
$$
Here $s_{n,d}$ is the number of near-$r$-regular rooted maps with $n$ edges and the degree of the root vertex equal to $d$. The numbers $d_{n,d}$ correspond to maps on the sphere with two distinct root vertices having total degree $d$. 

\begin{figure}[ht]
\centering
	\begin{subfigure}[b]{0.45\textwidth}
	\centering
    		\includegraphics[scale=2]{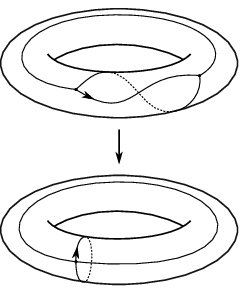}
		\caption{}
	\end{subfigure}
	\begin{subfigure}[b]{0.45\textwidth}
	\centering
    		\includegraphics[scale=2]{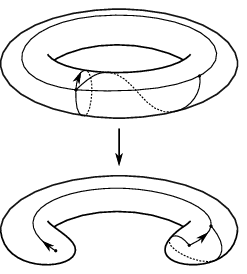}
		\caption{}
	\end{subfigure}	
	\caption{Toroidal maps}
\label{fig:linear}
\end{figure}

Indeed, let the root edge be a non-loop (see Figure \ref{fig:linear}(a)). After contracting such an edge we obtain some new rooted map with $n-1$ edges and the degree of its root vertex increased by two (summand $t_{n-1,d+r-2}$). Now suppose that the root edge is a loop. Two cases are possible: the root edge either encloses some region homeomorphic to a disc on the torus' surface or wraps the torus like a meridian or parallel (Figure \ref{fig:linear}(b)). 

In the first case contracting the root edge splits the torus into a new torus and a sphere, and the total degree of the root vertices of maps on these surfaces becomes equal to $d-2$. If the number of semi-edges incident to the root vertex and lying outside of the loop is $j$, then the number of semi-edges lying inside the loop is $d-j-2$. Similarly, if there are $i$ edges in total outside of the loop, then inside of it there are $n-i-1$ edges. Multiplying the number $t_{i,j}$ of ways to build a map on the torus with $i$ edges and $j$ semi-edges incident to the root vertex by the number $s_{n-1-i,d-j-2}$ of ways to build a map on the sphere with $n-1-i$ edges and a root vertex of degree $d-j-2$, we obtain the term $2\sum_{i=0}^{n-1}\sum_{j=0}^{d-2}s_{n-i-1,d-j-2}\cdot t_{i,j}$ in the right side of (\ref{eq:t_n_d}). The factor $2$ in this term arises from the fact that the sphere and the torus may be positioned in two ways with respect to the loop: the sphere outside and the torus inside or vice versa.

In the second case the root edge is a loop that wraps the torus (see Figure \ref{fig:linear}(b)). By contracting such edge we obtain a surface homeomorphic to a sphere with a map embedded into it. This map has $n-1$ edges, two distinct root vertices with a total degree $d-2$ and remaining vertices of degree $r$ (the term $d_{n-1,d-2}$ in (\ref{eq:t_n_d})). Positions of roots at the root vertices allow us to reconstruct the original toroidal map in a unique way.

A similar technique can be used to derive recurrence relations describing the number of maps on the projective plane and the Klein bottle. Let $p_{n, d}$ and $b_{n, d}$ be the numbers of near-$r$-regular maps on the projective plane and the Klein bottle, respectively. Arguments similar to those for maps on the torus allow us to write out the following recurrence relations for the numbers $p_{n, d}$ and $b_{n, d}$:
\begin{equation}
\label{eq:p_n_d}
p_{n,d}=p_{n-1,d+r-2}+(d-1) \cdot s_{n-1,d-2}+2\sum\limits_{i=0}^{n-1}\sum\limits_{j=0}^{d-2}s_{n-i-1,d-j-2}\cdot p_{i,j},
\end{equation}
\begin{equation}
\label{eq:bk_n_d}
b_{n,d}=b_{n-1,d+r-2} + (d-1) \cdot p_{n-1,d-2} + d_{n-1,d-2}+ 
\sum\limits_{i=0}^{n-1}\sum\limits_{j=0}^{d-2} \bigl(p_{n-i-1,d-j-2}\cdot p_{i,j} +2\cdot s_{n-i-1,d-j-2}\cdot b_{i,j}\bigr).
\end{equation}

To obtain a closed system of equations in addition to (\ref{eq:t_n_d}), (\ref{eq:p_n_d}) and (\ref{eq:bk_n_d}) we should also derive recurrence relations for $s_{n,d}$ and $d_{n,d}$. Using similar considerations one can show that the corresponding recurrence relations take the form
\begin{equation}
\label{eq:s_n_d}
s_{n,d}=s_{n-1,d+r-2}+\sum\limits_{i=0}^{n-1}\sum\limits_{j=0}^{d-2}s_{n-i-1,d-j-2}\cdot s_{i,j},
\end{equation}
\begin{equation}
\label{eq:d_n_d}
d_{n,d}=d_{n-1,d+r-2}-\sum\limits_{i=1}^{r-2}q^{(i)}_{n-1,d+r-2-i}+\dfrac{d(d-1)}{2}s_{n-1,d-2}+2\sum\limits_{i=0}^{n-1}\sum\limits_{j=0}^{d-2}s_{n-i-1,d-j-2}\cdot d_{i,j}.
\end{equation}
Here $q^{(i)}_{n,d}$ is the number of maps with $n$ edges, two distinct root vertices of the degrees $d$ and $i$, and with all the other vertices having the degree $r$, satisfying the following recurrence relation:
\begin{equation}
\label{eq:q_l_n_d}
q^{(i)}_{n,d}=q^{(i)}_{n-1,d+r-2}+i\cdot s_{n-1,d+i-2}+2\sum\limits_{j=0}^{n-1}\sum\limits_{k=0}^{d-2}q^{(i)}_{j,k}\cdot s_{n-j-1,d-k-2}.
\end{equation}

To derive analytic solutions for these recurrence relations we should rewrite them in a form of a system of equations for generating functions. We begin with the numbers $s_{n,d}$ and introduce the generating functions
$$
S_n(x)=\sum\limits_{d=0}^{2n}s_{n,d}\,x^d,\qquad\qquad S(x,t)=\sum\limits_{n=0}^{+\infty}S_n(x)\,t^n.
$$
In terms of generating functions, equation (\ref{eq:s_n_d}) can be rewritten as follows:
\begin{equation}
\label{eq:S_x_t_total}
S(x,t)-1=\dfrac{t}{x^{r-2}}\left[S(x,t)-1-\sum\limits_{i=1}^{r-2}x^i\,S_i(t)\right]+t\,x^2\,S^2(x,t),\qquad S_i(t):=\sum\limits_{j=0}^{\infty}s_{n,j}\,t^j.
\end{equation}
The possibility of deriving explicit expressions for the numbers $t_{n,d}$, $p_{n,d}$ and $b_{n,d}$ heavily depends on the possibility of obtaining an analytic solution of (\ref{eq:S_x_t_total}). This equation is quadratic in $S(x,t)$, but in addition to $S(x,t)$ it also contains $r-2$ functions $S_i(t)$ that should be determined. 

Our next steps depend on parity of $r$. For even $r$ from the handshaking lemma it follows that the root vertex must be of an even degree too. Consequently, $s_{n,j}=0$ for all odd $j$. This circumstance, in its turn, allows to reduce the number of unknown functions $S_i(t)$ by a factor of two. Moreover, for even $r$ it is convenient to rewrite the equations using a new variable $z=x^2$. Consequently we obtain the following equation for determining the generating function $S(z,t)$:
\begin{equation}
\label{eq:S_z_t_total}
S(z,t)-1=\dfrac{t}{z}\left[S(z,t)-1-\sum\limits_{i=1}^{r/2-1}z^i\,S_{2i}(t)\right]+t\,z\,S^2(z,t).
\end{equation}

Note that in one of his first papers \cite{Tutte_slicings} devoted to spherical map enumeration Tutte obtained an analytic expression for the numbers $s_{n,j}$ for even $r$. In the same paper Tutte noted that for odd $r$ he faced some difficulties. Further analysis conducted in \cite{Bender_1994}, \cite{Polyn_eq} showed the theoretical possibility of obtaining an analytic solution of (\ref{eq:S_z_t_total}) for all even values of $r$. However, in practice, the corresponding solution is known only for $r=4$, when the equation (\ref{eq:S_z_t_total}) contains in addition to $S(z,t)$ only one unknown function $S_2(t)$. As a consequence, for now some explicit forms for the numbers $t_{n,k}$, $p_{n,k}$ and $b_{n,k}$ can only be obtained for $4$-regular maps.

For odd values of $r$ the only case in which an analytic expression for (\ref{eq:S_x_t_total}) is known is the case of $r=3$. For larger $r$ the equation (\ref{eq:S_x_t_total}) has too many additional unknown functions $S_i(t)$, and even a theoretical possibility of deriving an explicit solution of (\ref{eq:S_x_t_total}) for such $r$ remains unclear \cite{Bender_1994}, \cite{Polyn_eq}.

In the next section we will use the technique described in paper \cite{Bender_1988} to obtain analytical solutions of this system for $r=4$. Then we will briefly describe the analogous results for $r=3$. 

\section{Enumeration of $4$-regular rooted maps} 

For $r=4$ the equation (\ref{eq:S_z_t_total}) can be rewritten as
$$
(z-t)\cdot S(z,t)-t\cdot z^2\cdot S^2(z,t)=(z-t)-t\cdot z\cdot S_2(t).
$$
For further reasoning it will be convenient to complete the square in the left hand side and rewrite this equation as
\begin{equation}
\label{eq:A2ztBzt}
A^2(z,t)=B(z,t).
\end{equation}
Here 
$$
A(z,t):=(z-t)-2t\,z^2\,S(z,t),\qquad \qquad B(z,t)=(z-t)^2-4tz^2(z-t)+4t^2z^3\,S_2(t).
$$
In these terms the recurrence relations (\ref{eq:t_n_d})--(\ref{eq:bk_n_d}) corresponding to the generating functions $T(z,t)$, $P(z,t)$, $B(z,t)$ take the form
\begin{equation}
\label{eq:T_z_t}
T(z,t)\cdot A(z,t)=t\,z\left[z\,D(z,t)-T_2(t)\right],
\end{equation}
\begin{equation}
\label{eq:P_zt}
P(z,t)\cdot A(z,t)=t\,z\left[z\,S(z,t)+2\,z^2\,\pd{S}{z}-P_2(t)\right],
\end{equation}
\begin{equation}
\label{eq:B_zt}
K(z,t)\cdot A(z,t)=t\,z\left[z\,P(z,t)\,(P(z,t)+1)+z\,D(z,t)+2\,z^2\,\pd{P}{z}-K_2(t)\right].
\end{equation}

To find the unknown function $S_2(t)$ we use a method suggested by Brown \cite{Brown_square_roots} (see also \cite{Bender_1988}, \cite{Polyn_eq}). Let $z=z_a(t)$ be the formal power series that zeroes out the value of $A(t,z_a(t),S(z_a(t),t))$. From the equality $A^2(z,t)=B(z,t)$ it follows that for $z=z_a(t)$ 
$$
B(z_a(t),t)=0,\qquad\qquad \left.\pd{B(z,t)}{z}\right|_{z=z_a(t)}=0.
$$
The two latter relations together with the equality $A(t,z_a(t),a(t))=0$, $a(t):=S(z_a(t),t)$, allow us to obtain an algebraic system of equations
\begin{equation}
\label{eq:tbz_a_system}
\begin{aligned}
&tz^2a^2-(z-t)a+(z-t)-tzS_2(t)=0,\\
&2tz^2a+t-z=0,\\
&2tza^2-a+1-tS_2(t)=0
\end{aligned}
\end{equation}
for determining the functions $z_a(t)$, $a(t)$, $S_2(t)$. It can be rewritten in the following parametric form:
\begin{equation}
\label{eq:tbz_a}
t=\dfrac{\sqrt{a-1}}{3a-2},\qquad S_2(t)=\sqrt{a-1}(2-a),\qquad z_a(t)=\dfrac{\sqrt{a-1}}{a}.
\end{equation} 
Then we can obtain an explicit expression for $S_2(t)$:
$$
S_2(t)=\dfrac{18\,t^2-1+(1-12\,t^2)^{3/2}}{54\,t^3}=t+2t^3+9t^5+54t^7+\ldots=
\sum\limits_{n=0}^{\infty}\sigma_{2n}^{(4)}\cdot t^{2n+1}.
$$
The numbers $\sigma_{2n}^{(4)}$ of $4$-regular maps with $n$ vertices and $2n$ edges of the sphere can be calculated by the formula
\begin{equation}
\label{eq:s_2n_4}
\sigma_{2n}^{(4)}=\dfrac{2\cdot 3^n\cdot(2n)!}{n!\cdot (n+2)!}
\end{equation}
(sequence A000168 on oeis.org). The function $S(z,t)$ is given by
$$
S(z,t)=\dfrac{z-t+\sqrt{(z-t)^2-4\,t\,z^2(z-t-z\,t\,S_2(t))}}{2\,t\,z^2}.
$$

From the formulas (\ref{eq:T_z_t})--(\ref{eq:B_zt}) we are essentially interested in finding explicit expressions for the functions $T_2(t)$, $P_2(t)$ and $B_2(t)$. They determine numbers of near-$4$-regular maps with the root of degree two and $n$ edges which are in bijection with rooted $4$-regular maps with $n-1$ edges. Since for $z=z_a(t)$ the function $A(z,t)$ in the left hand sides of (\ref{eq:T_z_t})--(\ref{eq:B_zt}) zeroes out, we obtain the following set of expressions:
\begin{equation}
\label{eq:T_2_t_eq}
T_2(t)=z_a(t)\cdot D(z_a(t),t),
\end{equation}
\begin{equation}
\label{eq:P_2_t_eq}
P_2(t)=a\,z_a(t)+2\,z_a^2(t)\left.\dfrac{\partial S(z,t)}{\partial z}\right|_{z_a(t)},
\end{equation}
\begin{equation}
\label{eq:K_2_t_eq}
K_2(t)=\left.\left(z\,P(z,t)\,(P(z,t)+1)+z\,D(z,t)+2\,z^2\,\pd{P}{z}\right)\right|_{z_a(t)}.
\end{equation}

To find an explicit form of $T_2(t)$ we need to know the value of $D(z,t)$ for $z=z_a(t)$. The equality (\ref{eq:d_n_d}) implies the following expression for $D(x,t)$ written in terms of the initial variables $(x,t)$:
$$
D(x,t)\cdot \left[x^2-t-2tx^4S(x,t)\right]=\dfrac{tx^4}{2}\pdd{}{x}{x}\left(x^2S(x,t)\right)-t\,x\,Q^{(1)}(x,t)-t\,x^2\,Q^{(2)}(x,t).
$$
Equations for the functions $Q^{(i)}(x,t)$, $i=1,2$ can be derived from the recurrence relations (\ref{eq:q_l_n_d}):
$$
Q^{(1)}(x,t)=\dfrac{t}{x^2}\left[Q^{(1)}(x,t)-x\,Q_1^{(1)}(t)\right]+t\,x\,S(x,t)+2\,t\,x^2\,S(x,t)\,Q^{(1)}(x,t),
$$
$$
Q^{(2)}(x,t)=\dfrac{t}{x^2}\left[Q^{(2)}(x,t)-x^2\,Q_2^{(2)}(t)\right]+2\,t\,S(x,t)+2\,t\,x^2\,S(x,t)\,Q^{(2)}(x,t).
$$
We rewrite these equations using instead of the variable $x$ and the functions $D(x,t)$, $S(x,t)$, $Q^{(1)}(x,t)$ and $Q^{(2)}(x,t)$ the variable $z=x^2$ and the functions $D(z,t)$, $A(z,t)$, $\tilde{Q}^{(1)}(z,t)=x(z)\,Q^{(1)}(x(z),t)$ and $Q^{(2)}(z,t)$. That results in the following set of equations:
\begin{equation}
\label{eq:D_z_t}
D(z,t)\cdot A(z,t)=-\dfrac{3}{2}\,A(z,t)+\dfrac{3}{2}\,z\,\pd{A(z,t)}{z}-z^2\,\pdd{A(z,t)}{z}{z}-\dfrac{3}{2}\,t-t\,\tilde{Q}^{(1)}(z,t)-t\,z\,Q^{(2)}(z,t),
\end{equation}
\begin{equation}
\label{eq:Q_1_z_t}
\tilde{Q}^{(1)}(z,t)\cdot A(z,t)=\dfrac{1}{2}\left(z-t-A(z,t)\right)-t\,z\,Q^{(1)}_1(t),
\end{equation}
\begin{equation}
\label{eq:Q_2_z_t}
\tilde{Q}^{(2)}(z,t)\cdot A(z,t)=t\,z\,\left[2\,S(z,t)-Q^{(2)}_2(t)\right].
\end{equation}

To determine the value of $D(z,t)$ for $z=z_a(t)$ we should take the derivative of (\ref{eq:D_z_t}) with respect to $z$ and find the limit as $z\to z_a(t)$:
\begin{equation}
\label{eq:D_z_t_z_a}
\left.D(z,t)\right|_{z=z_a(t)}=\dfrac{\left.\biggl[-\dfrac{1}{2}\,z\,\dfrac{\partial^2 A(z,t)}{\partial z^2}-
z^2\,\dfrac{\partial^3 A(z,t)}{\partial z^3}-
t\,\dfrac{\partial \tilde{Q}^{(1)}(z,t)}{\partial z}-t\,Q^{(2)}(z,t)-
t\,z\,\dfrac{\partial Q^{(2)}(z,t)}{\partial z}\biggr]\right|_{z=z_a(t)}}
{\left.\dfrac{\partial A(z,t)}{\partial z}\right|_{z=z_a(t)}}.
\end{equation}
Consequently, to obtain the value of $\left.D(z,t)\right|_{z=z_a(t)}$ we need to know the values of the functions $\tilde{Q}^{(1)}(z,t)$, $\tilde{Q}^{(2)}(z,t)$, their first derivatives, as well as the derivatives of $A(z,t)$ up to the third order inclusively at $z=z_a(t)$.

The derivatives of $A(z,t)$ at $z=z_a(t)$ can be obtained by sequential differentiation of the equality (\ref{eq:A2ztBzt}) with respect to $z$ taking into account that $A(z,t)=0$ when $S(z,t)=a$:
$$
\left.\pd{A}{z}\right|_{z_a(t)}=\dfrac{\sqrt{a\,(4-3a)}}{2-3a},
\quad
\left.\pdd{A}{z}{z}\right|_{z_a(t)}=\dfrac{4\,a^2\,\sqrt{a-1}}{(2-3a)\,\sqrt{a\,(4-3a)}},
\quad
\left.\dfrac{\partial^3A}{\partial z^3}\right|_{z_a(t)}=\dfrac{12\,a^3\,(a-1)}{(2-3a)\,(4-3a)\,\sqrt{a\,(4-3a)}}.
$$
To determine $\tilde{Q}^{(1)}(z,t)$ we need the value of $Q^{(1)}_1(t)$ from the right hand side of (\ref{eq:Q_1_z_t}). At $z=z_a(t)$ the left hand side of (\ref{eq:Q_1_z_t}) is equal to zero, so for $Q^{(1)}_1(t)$ the following equality holds:
$$
Q^{(1)}_1(t)=a\cdot z(a)=\sqrt{a-1}=\dfrac{1-\sqrt{1-12t^2}}{6t}=\sum\limits_{n=0}^{\infty}\rho^{(4)}_{2n}\cdot t^{2n+1}
=\sum\limits_{n=0}^{\infty}\dfrac{3^n}{n+1}\cdot\BCf{2n}{n}\cdot t^{2n+1}.
$$
Now take the derivative of the right hand side of (\ref{eq:Q_1_z_t}) with respect to $z$ and divide it by $\partial A/\partial z$ at $z=z_a(t)$:
$$
\left.\tilde{Q}^{(1)}(z,t)\right|_{z=z_a(t)}=-\dfrac{1}{2}+\dfrac{\dfrac{1}{2}-t\,Q^{(1)}_1(t)}{\left.\dfrac{\partial A}{\partial z}\right|_{z=z_a(t)}}=
-\dfrac{1}{2}+\dfrac{\sqrt{a}}{2\sqrt{4-3a}}.
$$
Finally, to determine the derivative of $\tilde{Q}^{(1)}(z,t)$ at $z=z_a(t)$ we differentiate the equation (\ref{eq:Q_1_z_t}) by $z$ one more time:
$$
\pdd{\tilde{Q}^{(1)}(z,t)}{z}{z}\,A(z,t)+2\,\pd{\tilde{Q}^{(1)}(z,t)}{z}\,\pd{A(z,t)}{z}+\tilde{Q}^{(1)}(z,t)\,\pdd{A(z,t)}{z}{z}=
-\dfrac{1}{2}\pdd{A(z,t)}{z}{z}.
$$
Thus, for the derivative of $\tilde{Q}^{(1)}(z,t)$ at $z=z_a(t)$ the following equality holds:
$$
\left.\dfrac{\partial\tilde{Q}^{(1)}(z,t)}{\partial z}\right|_{z=z_a(t)}=-\dfrac{1}{4}
\dfrac{\left.\dfrac{\partial^2 A}{\partial z^2}\right|_{z=z_a(t)}+2\left.\tilde{D}^{(1)}(z,t)\right|_{z=z_a(t)}}{\left.\dfrac{\partial A}{\partial z}\right|_{z=z_a(t)}}=
\dfrac{a\,\sqrt{a\,(a-1)\,(4-3a)}}{(4-3a)^2}.
$$
Using analogous considerations one can obtain the expressions for the function $Q^{(2)}(z,t)$ and for its derivative with respect to $z$ at $z=z_a(t)$:
$$
\left.Q^{(2)}\right|_{z=z_a(t)}=\sqrt{\dfrac{a}{a-1}}\,\left[\dfrac{2-a}{\sqrt{4-3a}}-1\right],
$$
$$
\left.\dfrac{\partial Q^{(2)}(z,t)}{\partial z}\right|_{z=z_a(t)}=\dfrac{a^2}{a-1}+\dfrac{a\,\sqrt{a}\,(a^2+2a-4)}{(a-1)\,(4-3a)\,\sqrt{4-3a}}.
$$
Substituting these expressions into the formula (\ref{eq:D_z_t_z_a}), we obtain an equality of the form
$$
\left.D(z,t)\right|_{z=z_a(t)}=\dfrac{a\,(a-1)}{(3a-4)^2},
$$
which allows us to use (\ref{eq:T_2_t_eq}) to find the value of $T_2(t)$:
$$
T_2(t)=z(a)\,\left.D(z,t)\right|_a=
\dfrac{t}{6}\,\left[\dfrac{1}{1-12\,t^2}-\dfrac{1}{\sqrt{1-12\,t^2}}\right]=\sum\limits_{n=0}^{\infty}\tau_{2n}^{(4)}\cdot t^{2n+1}=t^3+15t^5+198t^7+2511t^9+\ldots
$$
The numbers $\tau_{2n}^{(4)}$ enumerate $4$-regular maps with $n$ vertices and $2n$ edges on the torus. They can be written in an explicit form:
\begin{equation}
\label{eq:t_n_4}
\tau_{2n}^{(4)}=6^{n-1}\,\left[2^n-\dfrac{(2n-1)!!}{n!}\right],\qquad\qquad n=1,2,\ldots
\end{equation}

Using the relations derived above one can also apply the formula (\ref{eq:P_2_t_eq}) to obtain the following expression for the function $P_2(t)$ which defines the numbers $\pi_{2n}^{(4)}$ of $4$-regular maps on the projective plane:
$$
P_2(t)=\dfrac{1-\sqrt{a\,(4-3a)}}{\sqrt{a-1}}=t+\sum\limits_{n=1}^{\infty}\pi_{2n}^{(4)}\cdot t^{2n+1}=t+5\,t^3+38\,t^5+331\,t^7+3098\,t^9+\ldots
$$
Unfortunately, an explicit expression for $\pi_{2n}^{(4)}$ is quite cumbersome: using Lagrange-B\"urmann formula results in
$$
\pi_{2n}^{(4)}=\dfrac{3^{n}\cdot (2n)!}{(n+1)!\cdot n!}+\dfrac{3^{n+1}}{(2n+1)}
\sum\limits_{k=0}^{n-1}\BCf{2n+1}{k}\cdot2^{2k-2n}\cdot 
\left[\dfrac{1}{4}\cdot f(n+1-k)+\dfrac{4}{3}\cdot f(n-1-k)\right],
$$
where
$$
f(n)=\sum\limits_{i=0}^n\BCf{2i}{i}\cdot \BCf{2n-2i}{n-i}\cdot \dfrac{(-1)^i}{3^i}=
{}_2F_1\left(-n,\frac{1}{2};-n+\frac{1}{2};-\frac{1}{3}\right)\cdot\BCf{2n}{n}.
$$

From the formula (\ref{eq:K_2_t_eq}) it follows that to determine the numbers of $4$-regular maps on the Klein bottle we also need the values of $P(z,t)$ and its derivative with respect to $z$ at $z=z_a(t)$. Simple calculations show that
$$
\left.P(z,t)\right|_{z=z_a(t)}=-\dfrac{3\,\sqrt{a\,(4a-3)}+a-4}{2\,(4-3a)},
$$
$$
\left.\dfrac{\partial P(z,t)}{\partial z}\right|_{z=z_a(t)}=\dfrac{-a\sqrt{a-1}\,(2+2a+3\sqrt{a(4-3a)})}{(4-3a)^2}.
$$
Substituting these expressions into the right hand side of (\ref{eq:K_2_t_eq}), we obtain the generating function
$$
K_2(t)=\dfrac{2\sqrt{a-1}\,(4-a-3\sqrt{a(4-3a)}}{a\,(4-3a)^2}=\sum\limits_{n=0}^{\infty}\kappa_{2n}^{(4)}\cdot t^{2n+1}=4\,t^3+68\,t^5+964\,t^7+12836\,t^9+\ldots
$$
which determines the numbers $\kappa_{2n}^{(4)}$ of $4$-regular maps with $n$ vertices and $2n$ edges on the Klein bottle.

To enumerate unsensed $4$-regular maps on the torus we will also need some formulas for spherical maps with the root of degree $j$, two (the numbers $\hat{q}_{i,j}^{(2)}$) or three (the numbers $\hat{q}_{i,j}^{(3)}$) additional leaves and the remaining vertices of degree $4$. The corresponding recurrence relations take the form:
$$
\hat{q}_{i,j}^{(2)}=\hat{q}_{i-1,j+2}^{(2)}+q_{i-1,j-1}^{(1)}+
\sum\limits_{i=0}^{n-1}\sum\limits_{j=0}^{d-2}\left(2\,\hat{q}^{(2)}_{i,j}\cdot s_{n-i-1,d-j-2}+
q^{(1)}_{i,j}\cdot q^{(1)}_{n-i-1,d-j-2}\right),
$$
$$
\hat{q}_{i,j}^{(3)}=\hat{q}_{i-1,j+2}^{(3)}+\hat{q}_{i-1,j-1}^{(2)}+2
\sum\limits_{i=0}^{n-1}\sum\limits_{j=0}^{d-2}\left(\hat{q}^{(3)}_{i,j}\cdot s_{n-i-1,d-j-2}+
q^{(1)}_{i,j}\cdot \hat{q}^{(2)}_{n-i-1,d-j-2}\right),
$$
$$
\hat{q}_{i,j}^{(2)}=1\qquad\text{for $i=j=2$,}\qquad\qquad  \hat{q}_{i,j}^{(3)}=1\qquad\text{for $i=j=3$.}
$$
These recurrence relations can be transformed into the following equations for the generating functions $\hat{Q}^{(2)}(z,t)$ and $\hat{Q}^{(3)}(z,t)$:
$$
\hat{Q}^{(2)}(z,t)\cdot A(z,t)=z\,t\,\left[\tilde{D}^{(1)}(z,t)+[\tilde{D}^{(1)}(z,t)]^2-\hat{Q}_2^{(2)}(t)\right],
$$
$$
\hat{Q}^{(3)}(z,t)\cdot A(z,t)=t\,z^2\,\hat{Q}^{(2)}(z,t)+2\,t\,z^2\,\tilde{D}^{(1)}(z,t)\cdot \hat{Q}^{(2)}(z,t)-t\,z\,\hat{Q}_1^{(3)}(t).
$$
From the first equation one can find the following expression for the function $\hat{Q}^{(2)}$ at $a$:
$$
\left.\hat{Q}^{(2)}\right|_a=t\,z\,\dfrac{\left.\dfrac{\partial \tilde{D}^{(1)}}{\partial z}\right|_a\cdot (1+2\left.\tilde{D}^{(1)}\right|_a)}
{\left.\dfrac{\partial A}{\partial z}\right|_a}=\dfrac{(a-1)^{3/2}\,\sqrt{a}}{(4-3a)^{3/2}}.
$$
From the second equation an explicit expression for the function $\hat{Q}_1^{(3)}(t)$ follows:
$$
\hat{Q}_1^{(3)}(t)=\left.\hat{Q}^{(2)}\right|_a\cdot z\cdot\left(1+2\left.\tilde{D}^{(1)}\right|_a\right)=\dfrac{(a-1)^2}{(4-3a)^3}.
$$
Expressing $a$ through $t$ we finally obtain a formula for the numbers $\omega_n$ describing spherical maps with $n$ vertices of degree $4$ and four leaves, one of which is the root:
$$
\hat{Q}_1^{(3)}(t)=\dfrac{t^2}{6}\left[\dfrac{1}{(1-12t^2)^{3/2}}-\dfrac{1}{1-12t^2}\right]\quad\Longrightarrow\quad
\omega_n=6^{n-1}\left[\dfrac{(2n+1)!!}{n!}-2^n\right]=1,21,342,5049,\ldots
$$

\section{Enumeration of $3$-regular maps} 

In this section we briefly describe the solutions of analogous problems for $3$-regular maps on surfaces. Instead of the equation (\ref{eq:A2ztBzt}) we have
\begin{equation}
\label{eq:A2xtBxt}
A^2(x,t)=B(x,t)
\end{equation}
for 
$$
A(x,t):=(x-t)-2t\,x^3\,S(x,t),\qquad \qquad B(z,t)=(x-t)^2-4tx^3(x-t)+4t^2x^4\,S_1(t),
$$
and instead of (\ref{eq:tbz_a_system}) we have the system
\begin{equation}
\label{eq:tbx_a_system}
\begin{aligned}
&tx^3a^2-(x-t)a+(x-t)-txS_1(t)=0,\\
&2tx^3a+t-x=0,\\
&3tx^2a^2-a+1-tS_1(t)=0
\end{aligned}
\end{equation}
with the parametric solution
$$
x_a^3=\dfrac{a-1}{2a^2},\qquad t=\dfrac{a(1-a)}{2(1-2a)^3},\qquad 
S_1(t)=\dfrac{(2-a)(1-2a)^2}{a}.
$$
To determine an explicit form of the coefficients $s_n$ of the series expansion of $S_1(t)$ we introduce the parameter $\theta$ related to $t$ by
$$
\theta=-t^3\,r(\theta)\qquad\qquad\text{for}\qquad\qquad r(\theta)=\dfrac{1}{(1-2\theta)(1-4\theta)}
$$
and use the Lagrange-B\"urmann inversion formula. As a consequence we obtain that
$$
-[t^n]\theta=\dfrac{2^{2n}\,(3n)!!}{(n+1)!\,n!!},
$$
so the number $\sigma_{3n}^{(3)}$ of $3$-regular maps with $2n$ vertices and $3n$ edges can be calculated by the formula
$$
\sigma_{3n}^{(3)}=-\dfrac{2}{n+2}[t^n]\theta=\dfrac{2\cdot 4^n\cdot (3n)!!}{n!!\,(n+2)!}=1,4,32,336,\ldots
$$
(sequence $A002005$ on oeis.org).  

For $r=3$  instead of the equation (\ref{eq:T_2_t_eq}) we have the following system:
\begin{equation}
\label{eq:T_x_t_3}
T(x,t)\cdot A(x,t)=t\,x\,\left[x^2\,D(x,t)-T_1(t)\right],
\end{equation}
\begin{equation}
\label{eq:P_x_t_3}
P(x,t)\cdot A(x,t)=t\,x\,\left[x^2\,\pd{}{x}(x\,S(x,t))-P_1(t)\right],
\end{equation}
\begin{equation}
\label{eq:K_x_t_3}
K(x,t)\cdot A(x,t)=t\,x\,\left[x^2\,D(x,t)+x^2\,P(x,t)\,(P(x,t)+1)+x^3\pd{P}{x}-K_1(t)\right].
\end{equation}

Since the left hand side of (\ref{eq:T_x_t_3}) zeroes out at $x=x_a(t)$, so the right hand side does. Consequently, the generating function $T_1(t)$ for the numbers $t_n$ can be expressed through the function $a=a(t)$ by the formula
$$
T_1(t)=\left.x^2\,D(x,t)\right|_{x=x_a(t)}=\dfrac{t^5\,(2a-1)^5}{(1+2a-2a^2)^2}=
\dfrac{1}{(1-12\theta(1-2\theta))^2}\cdot\dfrac{1}{1-4\theta}.
$$
Applying the Lagrange-B\"urmann formula we conclude that the numbers $\tau_{3n}^{(3)}$ of $3$-regular maps having $2n$ vertices and $3n$ edges on the torus are equal to
$$
\tau_{3n}^{(3)}=[\lambda^n]\dfrac{1}{1-12\lambda(1-2\lambda)}\cdot\dfrac{1}{(1-4\lambda)^{n+2}(1-2\lambda)^{n+1}}=
$$
$$
=\dfrac{2^n}{(n+1)!}\sum\limits_{k=0}^{n}\dfrac{6^k}{(n-k)!}\,\sum\limits_{i=0}^{n-k}2^i\,\dfrac{(n+i+1)!}{i!}\,\dfrac{(2n-2k-i)!}{(n-k-i)!}=
$$
$$
=\dfrac{2^{2n}\,n!!}{(n+1)!}\sum\limits_{k=0}^n\dfrac{3^k\,(3n-2k+1)!!}{(n-k)!}=1,28,664,14912,326496,\ldots
$$

By analogous means from the formula (\ref{eq:P_x_t_3}) we obtain the following expression for the numbers $\pi_{3n}^{(3)}$ of maps having $2n$ vertices and $3n$ edges on the projective plane:
$$
P_1(t)=\left.\left[x^2\,S(x,t)+x^3\,\pd{S}{x}\right]\right|_{x=x_a(t)}=
\dfrac{1-\sqrt{1+2a-2a^2}}{2\,(2a-1)\,t}=
$$
$$
=t^2+\sum\limits_{n=1}^{\infty}\pi_{3n}^{(3)}\cdot t^{3n+2}=t^2+9\,t^5+118\,t^8+1773\,t^{11}+\ldots
$$
Applying the Lagrange-B\"urmann formula yields
$$
\pi_{3n}^{(3)}=-\dfrac{2^{2n+1}\cdot (3n)!!}{(n+1)!\cdot n!!}+\dfrac{3\cdot 2^{2n}}{(n+1)!!}\sum\limits_{k=0}^n\dfrac{3^k\cdot (2k-1)!!\cdot (3n-2k-1)!!}{2^k\cdot k!\cdot (n-k)!}.
$$

Finally, for the numbers $\kappa_{3n}^{(3)}$ describing $3$-regular maps with $2n$ vertices on the Klein bottle we obtain the generating function
$$
K_1(t)=\dfrac{3\,t^2\,(2a-1)^2\,(1-\sqrt{1+2a-2a^2})}{(1+2a-2a^2)^2}=\sum\limits_{n=1}^{\infty}\kappa_{3n}^{(3)}\cdot t^{3n+2}=6\,t^5+174\,t^8+4236\,t^{11}+\ldots,
$$
$$
\kappa_{3n}^{(3)}=
(n+1)\cdot (2\theta_n-\pi_{3n}^{(3)})+\dfrac{3\cdot 2^{2n}}{(n-2)!!}\sum\limits_{k=0}^n \dfrac{3^k\cdot (3n-2k-2)!!}{(n-k)!}.
$$

\section{The basic principles of enumerating $r$-regular sensed maps}

To enumerate $r$-regular maps up to orientation-preserving homeomorphisms we will use the technique developed in the article \cite{Mednykh_Nedela}. In that paper authors derived the following formula for the number $\tilde{\tau}_n$ of sensed maps with $n$ edges on the surface $S$ of a given genus $g$:
\begin{equation}
\label{eq:Mednykh_Ned_2006}
\tilde{\tau}_n=\dfrac{1}{2n}\sum_{\substack{L\mid 2n\\L\cdot m=2n}}\sum_{O\in {\rm Orb}(S_g/\Z_{L})}h_O(m) \cdot
\Epi_o(\pi_1(O), \Z_{L}).
\end{equation}
Here $h_O(m)$ is the number of rooted quotient maps with $m$ semi-edges on an orbifold $O$, corresponding to maps on $S$  with some predefined properties, $O$ runs over all orientation-preserving cyclic orbifolds ${\rm Orb}(S/\Z_{L})$ of the surface $S$ with period $L$, and $\Epi_o^{+}(\pi_1(O), \Z_{L})$ is the number of order-preserving epimorphisms from the fundamental group of the orbifold $O$ onto the cyclic group $\Z_{L}$. 

\begin{figure}[ht]
\centering
	\begin{subfigure}[b]{0.4\textwidth}
	\centering
    		\includegraphics[scale=1.3]{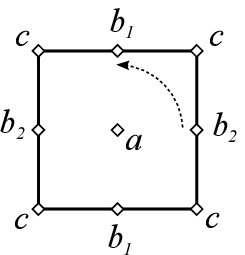}
		\caption{}
	\end{subfigure}	
	\begin{subfigure}[b]{0.4\textwidth}
	\centering
    		\includegraphics[scale=1.3]{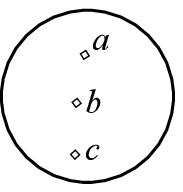}
		\caption{}
	\end{subfigure}	
	\caption{An automorphism of the torus and the corresponding orbifold}
\label{fig:orbifold}
\end{figure}

To illustrate these concepts consider the representation of a torus $T$ as a square with its opposite sides identified pairwise (Figure \ref{fig:orbifold} (a)). Rotation of this square by $90^\circ$ ($L=4$) splits the set of its points into two subsets, an infinite set of points in the general position and a finite set of singular points (see Figure \ref{fig:orbifold} (a)). Points in the general position are those that lie on an orbit of length $4$. Singular points are the remaining ones, and they necessary lie on an orbit of smaller length. For our example there are four singular points: $a$, $c$, $b_1$ and $b_2$. The former two of them are fixed, and the latter two are transformed into each other by the rotation by $90^\circ$. Identifying points of each orbit of this action, we obtain a sphere (Figure \ref{fig:orbifold} (b)). Critical points on the torus get transformed into {\em branch points} on the sphere (points $a, b, c$ on Figure \ref{fig:orbifold} (b)). From the topological point of view this construction can be viewed as a $4$-fold branched covering of the sphere $O$ by the torus $T$.  

For a general case of an orientation-preserving homeomorphism the corresponding manifold $O$ of genus $\mathfrak{g}$ with a finite number $r$ of branch points is called an {\em orbifold} with a signature
$$
O(\mathfrak{g} ,[m_1,\ldots,m_r]),\qquad 1<m_1\leq\ldots\leq m_r.
$$
Here $m_i$ are the {\em branch indices} of the corresponding branch points; each $m_i$ is equal to the period $L$ of the homeomorphism divided by the number of preimages of the corresponding branch point. For the example described above the branch points $a$ and $c$ have branch indices equal to $4$ and the branch index of $b$ is equal to $2$. Consequently, the signature of the corresponding orbifold takes the form
$$
O(0;[2,4,4])\equiv O(0;[2,4^2]).
$$
Rotation of the square by $270^\circ$ has an analogous description.

For the case of the torus there are five more periodic homeomorphisms that preserve its orientation and yield a sphere as an orbifold. The first of them corresponds to the rotation of the square by an angle of $180^\circ$. The remaining four are the rotations of a hexagon representing the torus by the angles of $60^\circ$, $120^\circ$, $240^\circ$ and $300^\circ$. The corresponding signatures have the form
$$
O(0;[2^4]),\qquad O(0;[2,3,6]), \qquad O(0;[3^3]), \qquad O(0;[3^3]), \qquad O(0;[2,3,6]).
$$
There also exists an infinite series of periodic homeomorphisms generated by two independent rotations of the torus, along its parallel and along its meridian. For a fixed period $L$ the number of such homeomorphisms is equal to $J_2(L)$, the second Jordan's function of $L$. As a consequence, the formula (\ref{eq:Mednykh_Ned_2006}) in the case of torus can be rewritten as
\begin{equation}
\label{eq:mednyh_main}
\tilde{\tau}_n^{(r)} = \dfrac{1}{2n}\biggl[\tilde{s}_{[2^4]}(n)+2\cdot\tilde{s}_{[2,4^2]}(n/2)+2\cdot\tilde{s}_{[3^3]}(2n/3)+
2\cdot\tilde{s}_{[2,3,6]}(n/3)+\sum_{L|n} J_2(L)\cdot \tau^{(r)}(n/L)\biggr].
\end{equation}
Here $\tilde{s}_{[m_1,\ldots,m_l]}(i)$ is the number of quotient maps with $i$ semi-edges on the sphere with branch points of indices $m_1,\ldots,m_l$, $\tau^{(r)}(i)$ is the number of rooted $r$-regular maps with $i$ edges on the torus. The latter numbers can be found using the formulas from the previous part of the present article. Hence the problem of enumerating sensed maps on the torus reduces to the problem of enumerating quotient maps on the corresponding orbifolds $O$. 

\begin{figure}[ht]
\centering
	\begin{subfigure}[b]{0.4\textwidth}
	\centering
    		\includegraphics[scale=1.3]{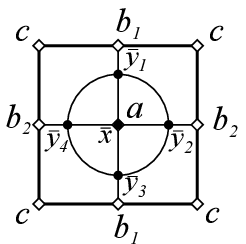}
		\caption{}
	\end{subfigure}
	\begin{subfigure}[b]{0.4\textwidth}
	\centering
    		\includegraphics[scale=1.3]{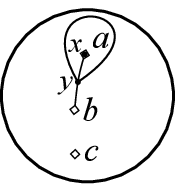}
		\caption{}
	\end{subfigure}
\caption{}
\label{fig:orbifold_map}
\end{figure}

As an example of such map consider a map $\M$ on the torus which is symmetric under the rotation of the square by $90^\circ$ (Figure \ref{fig:orbifold_map} (a)). As before, we identify the points lying on each orbit of the rotation, and instead of the map $\M$ obtain some {\em quotient map} $\mathfrak{M}$ on the orbifold  $O$ (Figure \ref{fig:orbifold_map} (b)). This quotient map would be a map on the sphere with the numbers of vertices, edges and faces equal to those of the original map $\M$ diveded by $4$ if the orbifold $O$ had no branch points and the surface $X$ had no corresponding critical points. The existence of such points makes the correspondence between these numbers more complicated.

Assume that a vertex $x$ of a quotient map $\mathfrak{M}$ coincides with some branch point of index $m_i$ on the orbifold $O$ (see vertex $x$ on Figure \ref{fig:orbifold}(b) which coincides with $a$). Then this vertex corresponds to $L/m_i$ vertices of the map $\M$ on the original surface $X$. The degree $d$ of $x$ in this case gets multiplied by $m_i$ on $X$ and becomes equal to $m_id$. For example, the vertex $x$ of degree $1$ of the quotient map $\mathfrak{M}$ shown on Figure \ref{fig:orbifold} (b) corresponds to a single vertex $\bar{x}$ of degree $4$ for the map $\M$ on Figure \ref{fig:orbifold}(a).

Now assume that a branch point of an orbifold $O$ falls into some face $f$ on the quotient map $\mathfrak{M}$ (see branch point $c$ on Figure \ref{fig:orbifold_map}(b)). This point will correspond to $L/m_i$ points on the torus $X$, $m_i$ being the corresponding branch index. The remaining points of $f$ are not branch points, so each of them corresponds to $L$ points on $X$. Hence, as in the case of a vertex, the degree of the face $f$ is multiplied by $m_i$ when this face is lifted to the manifold $X$. For example, the degree of the face that contains the branch point $c$ on Figure \ref{fig:orbifold_map}(b) is multiplied by $4$ on the torus.

One more peculiarity of quotient maps is the possibility of having {\em dangling semi-edges} which end not in vertices but in branch points of degree 2 (see branch point $b$ on Figure \ref{fig:orbifold_map}(b)). When lifted to the torus, any such edge gets transformed into $L/2$ edges of the map $\M$. These edges on the manifold $X$ contain critical points $b_i$ corresponding to the branch point $b$ (see critical points $b_1, b_2$ on Figure \ref{fig:orbifold_map} (a)). If an orbifold $O$ has no branch points of index $2$, then there are no dangling semi-edges for any quotient map $\mathfrak{M}$ on $O$. 

Thus to determine the functions $\tilde{s}_{[m_1,\ldots,m_l]}(i)$, where $i$ is the number of semi-edges, we need to iterate over all possible ways to distribute branch points of indices $m_1,\ldots,m_l$ over the faces, edges and vertices of a map, and for each of them determine the corresponding number of quotient maps on the sphere. As an example, in the next section we will perform these calculations for $4$-regular maps. 

\section{Enumeration of quotient maps of $4$-regular maps}

In the simplest case all branch points fall into faces of some quotient map $\mathfrak{M}$. This case can be realized for all four orbifolds:
$$
O(0;[2^4]),\qquad O(0;[2,4^2]), \qquad O(0;[3^3]), \qquad O(0;[2,3,6]).
$$
The corresponding quotient maps have $v'=v/L$ vertices of the degree $4$ that lift to $v$ vertices of degree $4$ of some map $\M$ on the torus. According to the Handshaking lemma, the number $n'$ of edges of a quotient map $\mathfrak{M}$ satisfies the equality
\begin{equation}
\label{eq:handshaking_lemma}
\sum\limits_{x\in V(\mathfrak{M})}\deg(x)=2n'
\end{equation}
and hence is equal to $2v'$. The number of faces $f'$ of a quotient map $\mathfrak{M}$ is defined by the Euler-Poincare formula
\begin{equation}
\label{eq:Euler_formula}
v'-n'+f'=2.
\end{equation}
Substituting $v'=v/L$, $n'=2v/L$, we find that $f'=2+v/L$.

Now we can count the number of such quotient maps $\mathfrak{M}$: it is the number  $\sigma(v')$ of $4$-regular maps on the sphere with $v'$ vertices multiplied by the number of ways to distribute branch points over $f'$ faces of each map $\mathfrak{M}$. This number should be calculated assuming that different branch points of the same index are indistinguishable. The results of the corresponding calculations are summarized in the following table:
$$
\begin{array}{ccccc}
\text{Orbifold}&L&v'&f'&\text{Number of quotient maps}\\[1,5 ex]
O(0;[2^4])&2&v/2&2+v/2&\sigma^{(4)}(v/2)\cdot\BCf{2+v/2}{4}\\[1,5 ex]
O(0;[2,4^2])&4&v/4&2+v/4&\sigma^{(4)}(v/4)\cdot(2+v/4)\cdot\BCf{1+v/4}{2}\\[1,5 ex]
O(0;[3^3])&3&v/3&2+v/3&\sigma^{(4)}(v/3)\cdot\BCf{2+v/3}{3}\\[1,5 ex]
O(0;[2,3,6])&6&v/6&2+v/6&\sigma^{(4)}(v/6)\cdot(2+v/6)\cdot(1+v/6)\cdot v/6\\[1,5 ex]
\end{array}
$$

A bit more complex is the case when some branch points coincide with vertices of a quotient map $\mathfrak{M}$ and the remaining branch points are located in its faces. For the orbifold $O(0;[3^3])$ this case can't be realized. For the orbifolds $O(0;[2,3,6])$, $O(0;[2,4^2])$ and $O(0;[2^4])$ it is possible that one or more branch points of index $2$ coincides with the vertices of degree $2$ of a quotient map $\mathfrak{M}$. Finally, for the case of the orbifold $O(0;[2,4^2])$ it is only possible that both of index $4$ branch points coincide with some vertices of $\mathfrak{M}$ of degree $1$. Indeed, if only one of branch points of index $4$ coincides with a vertex of degree $1$, the quotient map contains exactly one vertex of an odd degree which is impossible.

First we consider the case when only one branch point of degree $2$ coincides with a vertex $x$ of the degree $2$ of a quotient map $\mathfrak{M}$, whereas the remaining branch points are located in its faces. As we noted above, this case is possible for the orbifolds
$$
O(0;[2^4]),\qquad O(0;[2,4^2]), \qquad O(0;[2,3,6]).
$$
On the torus $T$ that covers the sphere $O$ the vertex $x$ corresponds to $L/2$ vertices of degree $4$, and each remaining vertex of $\mathfrak{M}$ corresponds to $L$ vertices of degree $4$. Consequently, the total number of vertices of a map $\mathfrak{M}$ is equal to
$$
v'=\dfrac{v-L/2}{L}+1=\dfrac{v+L/2}{L}.
$$
The number of its edges is calculated by (\ref{eq:handshaking_lemma}) as
$$
2n'=4\cdot\dfrac{v-L/2}{L}+2\cdot 1\qquad\Longrightarrow\qquad n'=\dfrac{2v}{L},
$$
and the number of its faces follows from the Euler-Poincare formula (\ref{eq:Euler_formula})
$$
f'=2-v'+n'=2+\dfrac{v-L/2}{L}.
$$
Note that we can replace a quotient map $\mathfrak{M}$ by some $4$-regular map $\M'$ on the sphere by replacing its vertex $x$ of the degree $2$ together with two edges incident to it by a single edge. The number of $4$-regular maps $\M'$ is equal to $\sigma^{(4)}((v-L/2)/L)$. To count quotient maps $\mathfrak{M}$ we note that any quotient map $\mathfrak{M}$ can be obtained from $\M'$ by choosing one of its edges in $n'-2$ ways and placing a vertex $x$ of degree $2$ in its middle and reassigning the root. Indeed, just placing the vertex $x$ yields a map that can't have its root at any of the edges incident to $x$. Since we should count maps with the root in an arbitrary position, we should use the double counting principle: dividing the number of  $4$-regular maps $\M'$ by $2(n'-2)$ and multiplying the result by $2n'$ results in exactly the number of $4$-regular maps with the vertex of degree $2$ and the root in an arbitrary position. To obtain the number of all quotient maps $\mathfrak{M}$ it remains to multiply the expression
$$
\dfrac{2n'}{2(n'-2)}\cdot \sigma^{(4)}((v-L/2)/L)
$$
by the number of ways to distribute branch points among $f'$ faces of a quotient map $\mathfrak{M}$. As a result, for quotient maps of the considered type we have the following values:
$$
\begin{array}{ccccc}
\text{Orbifold}&v'&f'&\text{Number of quotient maps}\\[1,5 ex]
(0,[2^4])&\dfrac{v-1}{2}&\dfrac{v+3}{2}&\sigma^{(4)}((v-1)/2)\cdot\BCf{(v+3)/2}{3}\cdot v \\[2 ex]
(0,[2,4^2])&\dfrac{v-2}{4}&\dfrac{v+6}{4}&\sigma^{(4)}((v-2)/4)\cdot \BCf{(v+6)/4}{2}\cdot \dfrac{v}{2}\\[2 ex]
(0,[2,3,6])&\dfrac{v-3}{6}&\dfrac{v+9}{6}&\sigma^{(4)}((v-3)/6)\cdot\dfrac{v+9}{6}\cdot \dfrac{v+3}{6}\cdot \dfrac{v}{3}\\[2 ex]
\end{array}
$$

Next we proceed with the case of at least one branch point of index $2$ coinciding with a vertex of a quotient map $\mathfrak{M}$. This is possible only for the orbifold $O(0;[2^4])$ and it remains to consider subcases of two, three and four index $2$ branch points coinciding with the vertices of $\mathfrak{M}$. Denote by $j$ the number of such branch points. The numbers of vertices, edges and faces of a quotient map can be calculated as
$$
v'=\dfrac{v+j}{2},\qquad n'=v,\qquad f'=2+\dfrac{v-j}{2}.
$$
Moreover, $j$ vertices have degree $2$ and the remaining $v'-j$ ones are of degree $4$. The problem of counting such maps again can be reduced to counting $4$-regular maps by replacing all degree $2$ vertices together with pairs of edges incident to them by single edges. The result of this operation is a $4$-regular map with $v'-j$ vertices and $n'-j$ edges. Hence to count quotient maps $\mathfrak{M}$ we start with the number $\sigma(v'-j)$ of $4$-regular maps, multiply it by the number $\BCCf{n'-j}{j}$ of ways to choose $j$ edges and place vertices on them. Then we use double counting to take into account that semi-edges incident to the new vertices can be roots too. The latter means multiplying the result by $2n'/(2n'-2j)$. Finally, multiplying the result by the number $\BCf{k'}{4-j}$ of ways to distribute the remaining $4-j$ branch points over $f'$ faces we obtain that the number of quotient maps $\mathfrak{M}$ is equal to
$$
\dfrac{2n'}{2n'-2j}\cdot\BCCf{n'-j}{j}\cdot\BCf{f'}{4-j}\cdot \sigma^{(4)}(v'-j)=
\BCf{v}{j}\cdot\BCf{(4+v-j)/2}{4-j}\cdot \sigma^{(4)}((v-j)/2).
$$

The final case in the group being considered is the one of both branch points of index $4$ on the orbifold $O(0;\{2,4^2\})$ coinciding with vertices of degree $1$ of a quotient map $\mathfrak{M}$, and the branch point of index $2$ coinciding with either its face or its edge. In the former sub-case the numbers $v'$ of vertices, $n'$ of edges and $f'$ of faces of a quotient map $\mathfrak{M}$ can be calculated as
$$
v'=\dfrac{v-2}{4}+2=\dfrac{v+6}{4},\qquad n'=\dfrac{v}{2},\qquad f'=2+n'-v'=\dfrac{v+2}{4}.
$$
To count such quotient maps we may utilize the fact that the number $\rho^{(4)}(v'-2)$ of rooted maps having $v'-2$ degree $4$ vertices and two leaves, one of which is a root, was already calculated in the third section. Taking all such maps, assigning new roots in them in $2n'$ ways instead of initial two ways and placing the branch point of degree $2$ into a face we obtain that the number of quotient maps $\mathfrak{M}$ is equal to
$$
\rho^{(4)}(v'-2)\cdot\dfrac{2n'}{2}\cdot f'=\dfrac{v(v+2)}{8}\cdot \rho^{(4)}((v-2)/4).
$$ 
Now assume that the branch point of degree $2$ coincides with a vertex too, namely with a vertex of degree $2$ of a map $\mathfrak{M}$. The numbers of vertices, edges and faces could be expressed as follows:
$$
v'=\dfrac{v-4}{4}+3,\qquad n'=\dfrac{v}{2},\qquad f'=2+n'-v'=\dfrac{v}{4}.
$$

Note that by replacing the vertex of degree $2$ together with the semi-edges incident to it we can obtain some map $\M'$ with $(v-4)/4+2$ vertices, namely with $(v-4)/4$ degree $4$ vertices and two leaves. Conversely, to transform the map $\M'$ back into $\mathfrak{M}$ we should choose an edge in $n'-1$ ways and place a new degree $2$ vertex $x$ on it. In the obtained map no semi-edge incident to $x$ could be its root, so to count the number of possible quotient maps $\mathfrak{M}$ we should start with the number $\rho^{(4)}(v'-3)$ of maps on the sphere with two leaves (one of which is the root) and $v'-3$ vertices of degree $4$, multiply it by the number $n'-1$ of ways to choose an edge in $\M'$ and adjust the result using double counting to take into account that the number of possible rootings should be $2n'$ instead of $2$. As a result we obtain the number of quotient maps in the form
$$
\rho^{(4)}(v'-3)\cdot (n'-1)\cdot \dfrac{2n'}{2}=\dfrac{v\,(v-2)}{4}\cdot \rho^{(4)}((v-4)/4).
$$

In the next group of cases at least one branch point of degree $2$ coincides with a dangling semi-edge of $\mathfrak{M}$. For the orbifold $O(0;[2,3,6])$ this case is impossible: it would result in a map of a single vertex of an odd degree. From the same consideration it follows that for the orbifold $O(0;\{2,4^2\})$ this sutuation is possible only if one out of two branch points of index $4$ coincides with a vertex and the other is located in a face. To count such quotient maps we consider maps $\M'$ on a sphere with $(v-1)/4$ vertices of degree $4$ and two leaves. The number $n'$ of edges and the number $f'$ of faces follow from the formulas (\ref{eq:handshaking_lemma}) and (\ref{eq:Euler_formula}): 
$$
2n'=4\dfrac{v-1}{4}+2=v+1,\qquad
f'=2+n'-v'=\dfrac{v+3}{4}.
$$
A quotient map $\mathfrak{M}$ is obtained from $\M'$ by removing one leaf together with a semi-edge incident to it. This procedure creates a dangling semi-edge in $\mathfrak{M}$. To count quotient maps we start with the number $\rho^{(4)}((v-1)/4)$ of spherical maps with two leaves, multiply it by the number $f'$ of ways to place a branch point into a face and use double counting because the quotient map may have its root in only $2n'-1$ positions. As a consequence, the number of such quotient maps can be calculated by the formula
$$
\rho^{(4)}((v-1)/4)\,\dfrac{v\,(v+3)}{4}.
$$

For the orbifold $O(0;\{2^4\})$ the fact that the number of odd-degree vertices must be even leaves only two possibilities: either two or four branch points of index $2$ fall into dangling semi-edges of a map $\mathfrak{M}$. The first case splits into three sub-cases depending on the placement of the remaining branch points: they are located in faces, in vertices, or in a face and a vertex of the quotient map $\mathfrak{M}$.

In the first sub-case we replace a quotient map $\mathfrak{M}$ with a combinatorially equivalent map $\M'$ with $v'$ vertices, $n'$ edges and $f'$ faces,
$$
v'=v/2+2,\qquad\qquad n'=v+1,\qquad\qquad f'=\dfrac{v+2}{2},
$$
which has $v/2$ vertices of degree $4$ and two vertices of degree $1$ which cannot be roots. To enumerate such maps we use the number $\rho^{(4)}(v/2)$ of maps with $v/2$ degree $4$ vertices, two leaves and a root in a leaf. The number of rootings of a quotient map should be recalculated since it's equal to $2n'-2=2v$ instead of $2$. We should also place two remaining branch points into $f'$ faces in $\BCf{f'}{2}$ ways. As a result, we obtain the number of quotient maps as
$$
\rho^{(4)}(v/2)\cdot v\cdot \BCf{(v+2)/2}{2}.
$$

In the second sub-case we considers maps $\M'$ which are equivalent to quotient maps $\mathfrak{M}$ and have the parameters
$$
v'=\dfrac{v-2}{2}+4,\qquad\qquad n'=v+1,\qquad\qquad f'=\dfrac{v}{2}.
$$
They have $v'-4$ vertices of degree $4$, two vertices of degree $2$, and two vertices of degree $1$. The latter ones couldn't be roots of a map. The number of such maps can be obtained from the number $\rho^{(4)}(v'-2)$ of maps with two leaves, $v'-2$ degree $4$ vertices and $n'-2$ edges by multiplying it by $\BCCf{n'-2}{2}$ ways to place two degree $2$ vertices on its edges and recalculating the number of ways to choose a root in such a map using double counting. Note that there are two possible rootings in maps enumerated by $\rho^{(4)}(v'-2)$ but the resulting maps have $2n'-2$ semi-edges which are incident to non-leaves and could all be roots. Consequently in this sub-case the number of quotient maps $\mathfrak{M}$ can be expressed as
$$
\rho^{(4)}((v-2)/2)\cdot v\cdot \BCCf{v-1}{2}=\rho^{(4)}((v-2)/2)\cdot v\cdot \BCf{v}{2}.
$$

In the third sub-case a quotient map $\mathfrak{M}$ is equivalent to a map $\M'$ with the parameters
$$
v'=\dfrac{v-1}{2}+3,\qquad\qquad n'=v+1,\qquad\qquad f'=\dfrac{v+1}{2},
$$
$(v-1)/2$ vertices of degree $4$, one vertex of degree $2$ and two non-root leaves. Such maps could be obtained by taking $\rho^{(4)}((v-1)/2)$ maps with $v'$ vertices of degree $4$, two leaves and $v$ edges, placing a degree $2$ vertex on some edge in $v$ ways, placing a new root in one of $2n'-2$ semi-edges, removing the existing root from one of two leaves, and labelling one of $f'$ faces with a branch point. As a consequence, the number of quotient maps $\mathfrak{M}$ is equal to
$$
\rho^{(4)}((v-1)/2)\cdot v^2\cdot \dfrac{v+1}{2}.
$$

The last case to consider is the one of the orbifold $O(0;\{2^4\})$ and all four branch points of degree $2$ coinciding with the dangling ends of semi-edges of $\mathfrak{M}$. Consider a map $\M'$ with the parameters
$$
v'=\dfrac{v}{2}+4,\qquad\qquad n'=v+2,\qquad\qquad f'=\dfrac{v}{2}
$$
having $v/2$ vertices of degree $4$ and four non-root leaves. All such maps could be obtained from $\omega(v/2)$ maps with a root in one of their four leaves and $v/2$ vertices of degree $4$ in the following way: in any given map place a new root in any of $2n'-4$ semi-edges incident to degree $4$ vertices and remove the old root which was incident to one of the four leaves. The double counting argument yields that the number of quotient maps $\mathfrak{M}$ in this case equals to
$$
\omega(v/2)\cdot \dfrac{v}{2}, \qquad\qquad \omega(v)=6^{v-1} \Biggl[ \frac{(2v+1)!!}{v!}-2^v \Biggr].
$$

Summarizing all these results, we may rewrite the formula (\ref{eq:mednyh_main}) for the case of unlabelled $4$-regular maps on the torus as
$$
\tilde{\tau}^{(4)}(v) = \dfrac{1}{4v}\Biggl[ \biggl[
\BCf{2+v/2}{4}\cdot\sigma^{(4)}\biggl(\frac{v}{2}\biggr)+ \frac{v}{2}\cdot \omega\biggl(\frac{v}{2}\biggr) +
  v\cdot \BCf{v}{2} \cdot\rho^{(4)}\biggl(\frac{v-2}{2}\biggr)+ \dfrac{v^2\,(v+1)}{2}\cdot \rho^{(4)}\biggl(\frac{v-1}{2}\biggr) +
$$
$$
+   v\cdot \BCf{(v+2)/2}{2}\cdot\rho^{(4)}\biggl(\frac{v}{2}\biggr) +v \cdot\BCf{(v+3)/2}{3}\cdot \sigma^{(4)}\biggl(\frac{v-1}{2}\biggr)+\BCf{v}{2}\cdot\BCf{(v+2)/2}{2}\cdot \sigma^{(4)}\biggl(\frac{v-2}{2}\biggr)+
$$
$$
+\dfrac{v+1}{2}\cdot\BCf{v}{3}\cdot \sigma^{(4)}\biggl(\frac{v-3}{2}\biggr)+\BCf{v}{4}\cdot \sigma^{(4)}\biggl(\frac{v-4}{2}\biggr) \biggl] + 2 \cdot \biggl[\frac{v+8}{4}\cdot\BCf{1+v/4}{2} \cdot\sigma^{(4)}\biggl(\frac{v}{4}\biggr)+\dfrac{v\,(v+3)}{4}\cdot \rho^{(4)}\biggl(\frac{v-1}{4}\biggr) +
$$
$$
+\dfrac{v}{2}\cdot\BCf{(v+6)/4}{2}\cdot   \sigma^{(4)}\biggl(\frac{v-2}{4}\biggr)+ \dfrac{v(v+2)}{8}\cdot \rho^{(4)}\biggl(\frac{v-2}{4}\biggr) + \dfrac{v\,(v-2)}{4}\cdot \rho^{(4)}\biggl(\frac{v-4}{4}\biggr) \biggr] + 
$$
$$
+ 2 \cdot \biggl [ \BCf{2+v/3}{3}\cdot \sigma^{(4)}\biggl(\frac{v}{3}\biggr)\biggr] + 2 \cdot \biggl[ \frac{v+12}{6}\cdot\frac{v+6}{6}\cdot \frac{v}{6}\cdot \sigma^{(4)}\biggl(\frac{v}{6}\biggr)+ \dfrac{v+9}{6}\cdot \dfrac{v+3}{6}\cdot \dfrac{v}{3}\cdot \sigma^{(4)}\biggl(\frac{v-3}{6}\biggr)\biggr] +
$$
$$
+\sum_{L|v} J_2(L)\cdot \tau^{(4)}\biggl(\frac{v}{L}\biggr)\biggr].
$$

\section{Enumeration of quotient maps of $r$-regular maps for other values of $r$}

In this section we briefly describe the corresponding results regarding enumeration of $r$-regular sensed maps on the torus for some other values of $r$. The case $r=3$ is very similar to the case $r=4$. Using arguments analogous to those given above, we obtain the following formula for the numbers $\tilde{\tau}^{(3)}(v)$ of sensed maps on the torus:
$$
\tilde{\tau}^{(3)}(v)=\dfrac{1}{3v}\Biggl[ \biggl[ \BCf{2+v/4}{4}\cdot \sigma^{(3)}\biggl(\frac{v}{2}\biggr) + \sum_{i=1}^4\dfrac{3v}{2i}\cdot \BCf{2+(v-2i)/4}{4-i}\cdot q^{(i)}\biggl(\frac{v}{2}+i\biggr) \biggr] + 
$$
$$
+ 2 \cdot \biggl[ \BCf{2+v/6}{3}\cdot \sigma^{(3)}\biggl(\frac{v}{3}\biggr) + \sum_{i=1}^3\dfrac{v}{i}\cdot \BCf{(v+12-4i)/6}{3-i}\cdot q^{(i)}\biggl(\frac{n+2i}{3}\biggr) \biggr] + 
$$

$$
+ 2 \cdot \biggl[ \dfrac{3v}{4}\cdot \BCf{(v+12)/8}{2}\cdot q^{(1)}\biggl(\frac{v}{4}+1\biggr)+\biggl(2+\frac{v}{8}\biggr)\cdot\BCf{1+v/8}{2}\cdot \sigma^{(3)}\biggl(\frac{v}{4}\biggr) \biggr] + 
$$

$$
+ 2 \cdot \biggl[ \biggl(2+\frac{v}{12}\biggr)\cdot \biggl(1+\frac{v}{12}\biggr)\cdot \frac{v}{12}\cdot \sigma^{(3)}\biggl(\frac{v}{6}\biggr) + v\cdot \BCf{(v+16)/12}{2}\cdot q^{(1)}\biggl(\frac{v+4}{6}\biggr) + 
$$
$$
+ \dfrac{v}{2}\cdot \dfrac{v+10}{12}\cdot q^{(2)}\biggl(\frac{v+10}{6}\biggr) + v\cdot \BCf{(v+18)/12}{2}\cdot q^{(1)}\biggl(\frac{v}{6}+1\biggr) \biggr] + \sum_{L|v} J_2(L)\cdot \tau^{(3)}\biggl(\frac{v}{L}\biggr)  \Biggr].
$$
Here the numbers $q^{(i)}(v)$ describe the number of maps with $v$ vertices, $i$ leaves and a root in one of these leaves and can be calculated from the following recurrence relation:
$$
q^{(i)}(v) = \frac{3v-2i-4}{i-1} \cdot q^{(i-1)}(v-2); \qquad q^{(2)}(2) = 1.
$$

For other values of $r$ instead of analytical formulas we have to use the recurrence relations obtained in the first sections of the article.  For a fixed $r>4$ any branch point that coincides with a vertex must have an index dividing $r$. Branch points that fall into dangling ends of semi edges must be of index $2$, and branch points that fall into vertices may have an arbitrary index. Taking these restrictions into account one can explicitly enumerate all the ways to distribute the branch points among vertices, semi-edges and faces of a quotient map. 

For a fixed distribution we would have to enumerate maps that have, in addition to degree $r$ vertices, some fixed amount of special vertices with degrees equal to divisors of $r$, as well as some amount of dangling semi-edges. To reduce such quotient maps to normal maps we should replace all dangling ends of semi-edges by leaves and use double counting to recalculate the number of possible rootings of a map. Finally, the number of maps on a sphere with some special vertices of degrees different from $r$ can be computed using Tutte's recursive approach of edge contraction. Note that contracting an edge which joins the root with a special vertex should be described by a separate summand for every possible degree of this special vertex. For the case of the root vertex being a loop, one would also have to consider all possible ways to distribute special vertices among two spherical maps that are obtained after the root is contracted.

\section*{Conclusion}

The results of calculations by the formulas obtained in this part of the article that define the numbers of rooted $r$-regular maps on the torus, on the projective plane and on the Klein bottle are given in the Tables \ref{table_last_torus}, \ref{table_last_pp} and \ref{table_last_KB}. The first terms of the corresponding sequences coincide with the results of explicit generation of the corresponding structures. Some of these sequences also match with known results obtained before by various authors. The results of computations of the numbers $\tilde{\tau}^{(r)}(v)$ of sensed $r$-regular maps on the torus are given in the Table \ref{table:table_last}. They were also verified using explicit generation.

The work was supported by grant 17-01-00212 from the Russian Foundation for Basic Research.

\begin{table}[h!]
\begin{center}
\footnotesize
\begin{tabular}{c|cccc}
\midrule
$ v $ &\phantom{00000}${\tau}^{(3)}(2v)$\phantom{00000}&\phantom{0000}${\tau}^{(4)}(v)$\phantom{0000}&\phantom{000000000000}${\tau}^{(5)}(2v)$\phantom{000000000000}&\phantom{000000000}${\tau}^{(6)}(v)$\phantom{000000000}\\ 
\midrule
1 & 1 & 1 & 120 & 10\\
2 & 28 & 15 & 125280 & 800\\
3 & 664 & 198 & 120800160 & 58000\\
4 & 14912 & 2511 & 113579366400 & 4080000\\
5 & 326496 & 31266 & 105549958379520 & 283100000\\
6 & 7048192 & 385398 & 97452182769223680 & 19496000000\\
7 & 150820608 & 4721004 & 89611995665911173120 & 1336380000000\\
8 & 3208396800 & 57590271 & 82178813933957614141440 & 91320000000000\\
9 & 67968706048 & 700465482 & 75217069050598359088496640 & 6226591000000000\\
10 & 1435486650368 & 8501284530 & 68747100051073934332046868480 & 423871680000000000\\
\midrule
\end{tabular}
\caption{Rooted $r$-regular maps on the torus}
\label{table_last_torus}
\end{center}
\end{table}

\begin{table}[h!]
\begin{center}
\footnotesize
\begin{tabular}{c|cccc}
\midrule
$ v $ &\phantom{00000}${\pi}^{(3)}(2v)$\phantom{00000}&\phantom{0000}${\pi}^{(4)}(v)$\phantom{0000}&\phantom{000000000000}${\pi}^{(5)}(2v)$\phantom{000000000000}&\phantom{000000000}${\pi}^{(6)}(v)$\phantom{000000000}\\ 
\midrule
1 & 9 & 5 & 215 & 22\\
2 & 118 & 38 & 106820 & 864\\
3 & 1773 & 331 & 65476730 & 40512\\
4 & 28650 & 3098 & 44355884860 & 2075860\\
5 & 484578 & 30330 & 31871222091735 & 112225776\\
6 & 8457708 & 306276 & 23809740820038860 & 6289396632\\
7 & 151054173 & 3163737 & 18286634336378438820 & 361699896960\\
8 & 2745685954 & 33252050 & 14338651143931504204140 & 21210328632420\\
9 & 50606020854 & 354312946 & 11425366917170617116755180 & 1262859239910000\\
10 & 943283037684 & 3817498004 & 9221856681066077433854516240 & 76114899842912520\\
\midrule
\end{tabular}
\caption{Rooted $r$-regular maps on the projective plane}
\label{table_last_pp}
\end{center}
\end{table}

\begin{table}[h!]
\begin{center}
\footnotesize
\begin{tabular}{c|cccc}
\midrule
$ v $ &\phantom{00000}${\kappa}^{(3)}(2v)$\phantom{00000}&\phantom{0000}${\kappa}^{(4)}(v)$\phantom{0000}&\phantom{000000000000}${\kappa}^{(5)}(2v)$\phantom{000000000000}&\phantom{000000000}${\kappa}^{(6)}(v)$\phantom{000000000}\\ 
\midrule
1 & 6 & 4 & 610 & 42	\\
2 & 174 & 68 & 713230 & 3846\\
3 & 4236 & 964 & 730427830 & 300048\\
4 & 97134 & 12836 & 714985017230 & 22171638\\
5 & 2163636 & 165784 & 684597649115160 & 1595739432\\
6 & 47394444 & 2103788 & 647152118916722050 & 113110095540\\
7 & 1027091736 & 26396416 & 606713944500089445300 & 7939173652032\\
8 & 22094309934 & 328621604 & 565552162701658630787310 & 553477015433958\\
9 & 472740763236 & 4068021916 & 524985730432063320579176680 & 38395171272416568\\
10 & 10074173087364 & 50142879128 & 485790080943651818443229561080 & 2653558455501023196\\
\midrule
\end{tabular}
\caption{Rooted $r$-regular maps on the Klein bottle}
\label{table_last_KB}
\end{center}
\end{table}

\begin{table}[h!]
\begin{center}
\footnotesize
\begin{tabular}{c|cccc}
\midrule
$ v $ &\phantom{00000}$\tilde{\tau}^{(3)}(2v)$\phantom{00000}&\phantom{0000}$\tilde{\tau}^{(4)}(v)$\phantom{0000}&\phantom{000000000000}$\tilde{\tau}^{(5)}(2v)$\phantom{000000000000}&\phantom{000000000}$\tilde{\tau}^{(6)}(v)$\phantom{000000000}\\ 
\midrule
1 & 1 & 1 & 15 & 3\\
2 & 5 & 4 & 6423 & 81\\
3 & 46 & 23 & 4031952 & 3313\\
4 & 669 & 185 & 2839677570 & 171282\\
5 & 11096 & 1647 & 2111005408320 & 9444158\\
6 & 196888 & 16455 & 1624203259187196 & 541659909\\
7 & 3596104 & 169734 & 1280171373413389056 & 31819176850\\
8 & 66867564 & 1805028 & 1027235174396893007472 & 1902508129720\\
9 & 1258801076 & 19472757 & 835745211680299639976976 & 115307287484560\\
10 & 23925376862 & 212603589 & 687471000510964612782875472 & 7064528615347192\\
\midrule
\end{tabular}
\caption{Sensed $r$-regular maps on the torus}
\label{table:table_last}
\end{center}
\end{table}

\end{document}